# Multi-objective Geometric Programming Problem with Weighted Mean Method

A. K. Ojha
School of Basic Sciences
IIT Bhubaneswar
Orissa, Pin-751013, India
.

K. K. Biswal
Department of Mathematics
C.T.T.C, Bhubaneswar
Orissa, India-751024
.

*Abstract*-**Geometric programming is an important class of optimization problems that enable practitioners to model a large variety of real-world applications, mostly in the field of engineering design. In many real life optimization problem multi-objective programming plays a vital role in socio economical and industrial optimizing problems. In this paper we have discussed the basic concepts and principle of multiple objective optimization problems and developed geometric programming (GP) technique to solve this optimization problem using weighted method to obtain the non-inferior solutions.**

*Keywords- Geometric programming, optimization, weighted method, duality theory, non-inferior solutions.*

## 1 INTRODUCTION

Geometric Programming(GP) problems have wide range of application in production planning, location, distribution, risk managements, chemical process designs and other engineer design situations. GP problem is an excellent method when decision variable interact in a non-linear, especially, in an exponential fashion. Most of these GP applications are posynomial type with zero or few degrees of difficulty. GP problem whose parameters, except for exponents, are all positive are called posynomial problems, where as GP problems with some negative parameters are refereed to as signomial problems. The degree of difficulty is defined as the number of terms minus the number of variables minus one, and is equal to the dimension of the dual problem. When the degree of difficulty is zero, the problem can be solved analytically. For such posynomial problems, GP techniques find global optimal solutions.

If the degree of difficulty is positive, then the dual feasible region must be searched to maximize the dual objective, while if the degree of difficulty is negative, the dual constraints may be inconsistent. For detail discussions of various algorithms and computational aspects for both posynomial and signomial GP refers to Beightler et al.[2], Duffin et al.[6], Ecker[7] and Phillips et al.[13]. From early 1960 a lot of research work have been done on these GP problems [3, 5, 9, 11, 14, 16, 17, 18]. Mainly, we use GP technique to solve some optimal engineering design problems[1] where we minimize cost and/or weight, maximize volume and/or efficiency etc. Generally, an engineering design problem has multiple objective functions. In this case, it is not suitable to use any single-objective programming to find an optimal compromise solution. Biswal[4] has studied the optimal compromise solution of multi-objective programming problem by using fuzzy programming technique [20, 21]. In a recent paper Islam and Ray[8] find the Pareto optimal solution by considering a multi-objective entropy transportation problem with an additional restriction with generalized fuzzy number cost.

In this paper we have developed the method to find the compromise optimal solution of certain multi-objective geometric programming problem by using weighting method. First of all the multiple objective functions transformed to a single objective by considering it as the linear combination of the multiple objectives along with suitable constants called weights. By changing the weights the most compromise optimal solution has been arrived by using GP techniques.

The organization of the paper is as follows: Following the introduction, Formulation of multi-objective GP and corresponding weighting method have been discussed in Section 2 and 3. The duality theory has been discussed at Section 4 to find the optimal value of the objective function and the illustrative examples have been incorporated at Section 5 to understand the problem. Finally, at Section 6 some conclusions are drawn from the discussion.

## 2 FORMULATION OF MULTI-OBJECTIVE GEOMETRIC PROGRAMMING

A multi-objective geometric programming problem can be defined as:

Find $x = (x_1, x_2, ...., x_n)^T$

so as to

$$\min : g_{k0}(x) = \sum_{t=1}^{T_{k0}} C_{k0t} \prod_{j=1}^{n} x_j^{a_{k0tj}}, \quad k = 1, 2, ..., p \quad (2.1)$$

subject to

$$g_i(x) = \sum_{t=1}^{T_i} C_{it} \prod_{j=1}^{n} x_j^{d_{itj}} \leq 1, \quad i = 1, 2, ..., m \quad (2.2)$$

$$x_j > 0, \quad j = 1, 2, ..., n \quad (2.3)$$

where $C_{k0t}$ for all k and t are positive real numbers and $d_{itj}$ and $a_{k0tj}$ are real numbers for all i, k, t, j.

$T_{k0}$ = number of terms present in the $k^{th}$ objective function.
$T_i$ = number of terms present in the $i^{th}$ constraint.
In the above multi-objective geometric program there are p number of minimization type objective function, m number of inequality type constraints and n number of strictly positive decision variables.

## 3 WEIGHTING METHOD OF MULTI-OBJECTIVE FUNCTIONS:

The weighting method is the simplest multi-objective optimization which has been widely applied to find the non-inferior optimal solution of multi-objective function within the convex objective space.

If $f_1(x)$, $f_2(x)$,... ,$f_n(x)$ are n objective functions for any vector $x = (x_1, x_2, ...., x_n)^T$ then we can define weighting method for their optimal solution as defined below:

Let $W = \left\{ w : w \in R^n, w_i > 0, \sum_{j=1}^{n} w_j = 1 \right\}$ to be the set of non-negative weights. The weighted objective function for the multiple objective function defined above can be defined as P(w) where

$$P(w) = \min_{x \in X} \sum_{j=1}^{n} w_j f_j(x) \quad (3.1)$$

It must be made clear, however, that if the objective space of the original problem is non-convex, then the weighting method may not be capable of generating the efficient solutions on the non-convex part of the efficient frontier. It must also be noted that the optimal solution of a weighting problem should not be used as the best compromise solution, if the weights do not reflect the decision maker's preferences or if the decision maker does not accept the assumption of a linear utility function. For more detail about the weighted method refer[10].

Based on the importance of the p number of objective functions defined in (2.1) the weights $w_1, w_2, ..., w_p$ are assigned to define a new min type objective function Z(w) which can be defined as

$$\min : Z(x) = \sum_{k=1}^{p} w_k g_{k0}(x)$$

$$= \sum_{k=1}^{p} w_k \left( \sum_{t=1}^{T_{k0}} C_{k0t} \prod_{j=1}^{n} x_j^{a_{k0tj}} \right) \quad (3.2)$$

$$= \sum_{k=1}^{p} \sum_{t=1}^{T_{k0}} w_k C_{k0t} \prod_{j=1}^{n} x_j^{a_{k0tj}}$$

subject to $\sum_{t=1}^{T_i} C_{it} \prod_{j=1}^{n} x_j^{d_{itj}} \leq 1$, $i = 1, 2, ..., m$ (3.3)

$x_j > 0$, $j = 1, 2, ..., n$

where $\sum_{k=1}^{p} w_k = 1$, $w_k > 0$, $k = 1, 2, ..., p$ (3.4)

## 4. DUAL FORM OF GPP

The model given by (3.2), (3.3) and (3.4) is a conventional geometric programming problem and it can be solved directly by using primal based algorithm for non linear primal problem or dual programming [12]. Methods due to Rajgopal and Bricker [15], Beightler and Phillips[1] and Duffin et al.[6] projected in their analysis that the dual problem has the desirable features of being linearly constrained and having an objective function with structural properties with suitable solution.

According to Duffin et al.[6] the model given by (3.3) can be transformed to the corresponding dual geometric program as:

$$\max_{w} \prod_{t=1}^{T_0} \left( \frac{w_k C_{k0t}}{w_{0t}} \right)^{w_{0t}} \prod_{i=1}^{m} \prod_{t=1}^{T_i} \left( \frac{w_{i0} C_{it}}{w_{it}} \right)^{w_{it}} \prod_{t=1}^{T_i} \lambda(w_{it})^{\lambda(w_{it})}$$
(4.1)

subject to $\sum_{t=1}^{T_0} w_{0t} = 1$

$\sum_{i=1}^{m} \sum_{t=1}^{T_i} a_{itj} w_{it} + \sum_{i=1}^{m} \sum_{t=1}^{T_i} d_{itj} w_{it} = 0$, j=1,2,...,n

$w_{it} \geq 0 \quad \forall t, j$

$\sum_{k=1}^{p} w_k = 1$, $w_k > 0$, $k = 1, 2, ..., p$

Since it is a usual dual problem then it can be solved using method relating to the dual theory.

## 5 NUMERICAL EXAMPLES

For illustration we consider the following examples.
Example: 1
Find $x_1, x_2, x_3, x_4$ so as to

$\min : g_{10}(x) = 4x_1 + 10x_2 + 4x_3 + 2x_4$ (5.1)

$\max : g_{20}(x) = x_1 x_2 x_3$ (5.2)

subject to $\dfrac{x_1^2}{x_4^2} + \dfrac{x_2^2}{x_4^2} \leq 1$ (5.3)

$\dfrac{100}{x_1 x_2 x_3} \leq 1$ (5.4)

$x_1, x_2, x_3, x_4 > 0$

Now the problem can be rewritten as

$\min : g_{10}(x) = 4x_1 + 10x_2 + 4x_3 + 2x_4$ (5.5)

$\min : g'_{20}(x) = x_1^{-1} x_2^{-1} x_3^{-1}$ (5.6)

Subject to $x_1^2 x_4^{-2} + x_2^2 x_4^{-2} \leq 1$ (5.7)

$100 x_1^{-1} x_2^{-1} x_3^{-1} \leq 1$ (5.8)

$x_1, x_2, x_3, x_4 > 0$ (5.9)

Introducing weights for the above objective functions a new objective function is formulated as

$$Z(x) = w_1(4x_1 + 10x_2 + 4x_3 + 2x_4) + w_2(x_1^{-1}x_2^{-1}x_3^{-1})$$
(5.10)

subject to $x_1^2 x_4^{-2} + x_2^2 x_4^{-2} \leq 1$ (5.11)

$$100 x_1^{-1} x_2^{-1} x_3^{-1} \leq 1$$ (5.12)

$$x_1, x_2, x_3, x_4 > 0$$

where $w_1 + w_2 = 1, \quad w_1, w_2 > 0$ (5.13)

This problem is having degree of difficulty 3. The problem is solved via the dual programming [6].
The corresponding dual program is:

$$\max_w : V(w) = \left(\frac{4w_1}{w_{01}}\right)^{w_{01}} \left(\frac{10w_1}{w_{02}}\right)^{w_{02}} \left(\frac{4w_1}{w_{03}}\right)^{w_{03}} \left(\frac{2w_1}{w_{04}}\right)^{w_{04}}$$

$$\left(\frac{w_2}{w_{05}}\right)^{w_{05}} \left(\frac{1}{w_{11}}\right)^{w_{11}} \left(\frac{1}{w_{12}}\right)^{w_{12}} (w_{11}+w_{12})^{(w_{11}+w_{12})} 100^{w_{21}}$$
(5.14)

subject to $w_{01} + w_{02} + w_{03} + w_{04} + w_{05} = 1$

$$w_{01} - w_{05} + 2w_{11} - w_{21} = 0$$

$$w_{02} - w_{05} - 2w_{12} - w_{21} = 0$$

$$w_{03} - w_{05} - w_{21} = 0$$

$$w_{04} - 2w_{11} - 2w_{12} = 0$$

$$w_1 + w_2 = 1$$

$$w_{01}, w_{02}, w_{03}, w_{04}, w_{05}, w_{11}, w_{12}, w_{21} \geq 0$$

$$w_1, w_2 > 0$$

By considering different values of $w_1$ and $w_2$ the dual variables, corresponding maximum value of dual objective are given in the following Table.

Table-1
Dual solution:

| $w_1$ | $w_2$ | $w_{01}$ | $w_{02}$ | $w_{03}$ | $w_{04}$ |
|---|---|---|---|---|---|
| 0.1 | 0.9 | 0.2308894 | 0.3045667 | 0.3329927 | 0.1305293 |
| 0.2 | 0.8 | 0.2310206 | 0.3044397 | 0.3331819 | 0.1306035 |
| 0.3 | 0.7 | 0.2310643 | 0.3047974 | 0.3332450 | 0.1306282 |
| 0.4 | 0.6 | 0.2310862 | 0.3048263 | 0.3332765 | 0.1306406 |
| 0.5 | 0.5 | 0.2310993 | 0.3048436 | 0.3332955 | 0.1306480 |

| $w_{05}$ | $w_{11}$ | $w_{12}$ | $w_{21}$ | Z |
|---|---|---|---|---|
| 0.0010217 | 0.051051 | 0.014213 | 0.3319702 | 8.80776 |
| 0.0004543 | 0.051080 | 0.014221 | 0.3319701 | 17.60555 |
| 0.000265092 | 0.05109 | 0.0142237 | 0.3329799 | 26.40333 |
| 0.000170424 | 0.051095 | 0. 014225 | 0.3331061 | 35.201110. |
| 0.000113614 | 0.051098 | 0.0142259 | 0.3331818 | 43.99888 |

Using primal dual relationship the corresponding primal solution are given in the following Table.

Table-2
Primal solution:

| $w_1$ | $w_1$ | $x_1$ | $x_2$ | $x_3$ | $x_4$ | Z |
|---|---|---|---|---|---|---|
| 0.1 | 0.9 | 5.084055 | 2.682555 | 7.332315 | 5.748367 | 8.80776 |
| 0.2 | 0.8 | 5.084055 | 2.682555 | 7.332315 | 5.748367 | 17.60555 |
| 0.3 | 0.7 | 5.084055 | 2.682555 | 7.332315 | 5.748367 | 26.40333 |
| 0.4 | 0.6 | 5.084055 | 2.682555 | 7.332315 | 5.748367 | 35.20111 |
| 0.5 | 0.5 | 5.084055 | 2.682555 | 7.332315 | 5.748367 | 43.99888 |

where the minimum values of $g_{10} = 87.98776$ and $g_{20} = 0.01$

Example: 2
Find $x_1, x_2, x_3, x_4$ in order to

$$\min : f_1(x) = x_1^{-1} x_2^{-1} x_3^{-1}$$ (5.15)

$$\min : f_2(x) = x_1^{-1} x_2^{-3} x_3^{-5} + x_1^{-1} x_2^{-1}$$ (5.16)

subject to $x_1 x_2 x_3^2 + x_2 x_3 \leq 6$ (5.17)

$$x_1 x_3 \leq 1$$ (5.18)

$$x_1, x_2, x_3 > 0$$ (5.19)

Using the weights the above objective function can be reduced to the new objective function as:

$$Z(x) = w_1(x_1^{-1} x_2^{-1} x_3^{-1}) + w_2(x_1^{-1} x_2^{-3} x_3^{-5} + x_1^{-1} x_2^{-1})$$
(5.20)

$$x_1 x_2 x_3^2 + x_2 x_3 \leq 6$$ (5.21)

$$x_1 x_3 \leq 1$$ (5.22)

$$x_1, x_2, x_3 > 0$$ (5.23)

where $w_1 + w_2 = 1 \quad w_1, w_2 > 0$ (5.24)

In this problem the degree of difficulty is 2 and it can be solved by using duality theory as given by

$$\max_w : V(w) = \left(\frac{w_1}{w_{01}}\right)^{w_{01}} \left(\frac{w_2}{w_{02}}\right)^{w_{02}} \left(\frac{w_2}{w_{03}}\right)^{w_{03}}$$

$$\left(\frac{1}{6w_{11}}\right)^{w_{11}} \left(\frac{1}{6w_{12}}\right)^{w_{12}} (w_{11}+w_{12})^{(w_{11}+w_{12})}$$
(5.25)

subject to $w_{01} + w_{02} + w_{03} = 1$

$$-w_{01} - w_{02} - w_{03} + w_{11} + w_{21} = 0$$

$$-w_{01} - 3w_{02} - w_{03} + w_{11} + w_{12} = 0$$

$$-2w_{01} - 5w_{02} + 2w_{11} + w_{12} + w_{21} = 0$$

$$w_{01}, w_{02}, w_{03}, w_{11}, w_{12}, w_{21} \geq 0$$

$$w_1, w_2 > 0$$

For different values of $w_1$, $w_2$ the dual variables and the corresponding maximum values of dual objective is obtained as given in the Table

Table-3
Dual solution:

| $w_1$ | $w_2$ | $w_{01}$ | $w_{02}$ | $w_{03}$ |
|---|---|---|---|---|
| 0.1 | 0.9 | 0.2085711 | 0.5276192 | 0.2638096 |
| 0.2 | 0.8 | 0.3640122 | 0.4239919 | 0.2119959 |
| 0.3 | 0.7 | 0.4952513 | 0.3364992 | 0.1682496 |
| 0.4 | 0.6 | 0.604162 | 0.2638920 | 0.1319460 |
| 0.5 | 0.5 | 0.6959958 | 0.2026694 | 0.1013347 |

| $w_{11}$ | $w_{12}$ | $w_{21}$ | Z |
|---|---|---|---|
| 1.00 | 1.055235 | 0 | 0.1642316 |
| 0.9239919 | 0.9239918 | 0.076008 | 0.1831441 |
| 0.8364992 | 0.8364992 | 0.1635008 | 0.2019177 |
| 0.7638920 | 0.7638920 | 0.2361080 | 0.2206914 |
| 0.7026695 | 0.7026694 | 0.2973305 | 0.2394650 |

The corresponding primal solution is given in the following Table:

Table-4
Primal Solution

| $w_1$ | $w_1$ | $x_1$ | $x_2$ | $x_3$ | Z |
|---|---|---|---|---|---|
| 0.1 | 0.9 | 2.527860 | 8.217575 | 0.3748833 | 0.1642316 |
| 0.2 | 0.8 | 2.620746 | 7.862237 | 0.3815708 | 0.1831441 |
| 0.3 | 0.7 | 2.620745 | 7.862236 | 0.3815709 | 0.2019177 |
| 0.4 | 0.6 | 2.620747 | 7.862240 | 0.3815707 | 0.2206914 |
| 0.5 | 0.5 | 2.620747 | 7.862242 | 0.3815705 | 0.2394650 |

Ideal solution of the two objective functions are given below:
$f_1 = 0.33333$; $x_1 = 236.9322$; $x_2 = 710.7964$; $x_3 = 0.0042206$ and $f_2 = 0.1421595$; $x_1 = 2.148558$; $x_2 = 9.82199$; $x_3 = 0.3490711$.

6 CONCLUSIONS

By using weighted method we can solve a multi-objective GPP as a vector minimum problem. A vector-maximum problem can be transformed as a vector minimization problem. If any of the objective function and/or constraint does not satisfy the property of a posynomial after the transformation, then we use any of the general purpose non-linear programming algorithm to solve the problem. We can also use this technique to solve a multi-objective signomial geometric programming problem. However, if a GPP has either a higher degree of difficulty or a negative degree of difficulty, then we can use any of the general purpose non-linear programming algorithm instead of a GP algorithm.

AUTHORS PROFILE

**Dr.A.K.Ojha:** Dr.A.K.Ojha received a Ph.D(mathematics) from Utkal University in 1997. Currently he is an Asst.Prof. in Mathematics at I.I.T. Bhubaneswar, India. He is performing research in Nural Network, Geometric Programming, Genetical Algorithem, and Particle Swarm Optimization. He has served more than 27 years in different Govt. colleges in the state of Orissa. He has published 22 research papers in different journals and 7 books for degree students such as: Fortran 77 Programming, A text book of modern algebra, Fundamentals of Numerical Analysis etc.

**K.K.Biswal:** Mr.K.K.Biswal received M.Sc.(Mathematics) from Utkal University in 1996. Currently he is a lecturer in Mathematics at CTTC, Bhubaneswar, India. He is performing research works in Geometric Programming. He is served more than 7 years in different colleges in the state of Orissa. He has published 2 research papers in different journals.